\documentclass{amsart}[14pt]

\usepackage{a4wide}
\usepackage{latexsym}
\usepackage{amsfonts}
\usepackage{amsmath}
\usepackage{euscript}
\usepackage[all]{xy}
\usepackage{amscd}
\usepackage{amsthm}
\usepackage{amssymb}
\usepackage{verbatim}

\usepackage[T1]{CJKutf8}

\newtheorem{theorem}{Theorem}[subsection]
\newtheorem{proposition}[theorem]{Proposition}
\newtheorem{lemma}[theorem]{Lemma}
\newtheorem{corollary}[theorem]{Corollary}

\theoremstyle{definition}
\newtheorem{definition}[theorem]{Definition}
\newtheorem{example}[theorem]{Example}

\theoremstyle{remark}
\newtheorem{remark}[theorem]{Remark}

\numberwithin{equation}{section}

\def\Ob{\mathop{\rm Ob}}
\def\lim{\mathop{\varprojlim}}

\def\Hom{\mathop{\rm Hom}}
\def\Ext{{\mathop{\rm Ext}}}
\def\k{\underline{k}}
\def\End{{\mathop{\rm End}}}
\def\Aut{{\mathop{\rm Aut}}}
\def\Res{{\mathop{\rm Res}}}
\def\res{{\mathop{\rm res}}}
\def\Mor{\mathop{\rm Mor}}

\def\Id{\mathop{\rm Id}}

\def\C{\mathcal{C}}
\def\c{\mathbb{C}}
\def\D{\mathcal{D}}
\def\d{\mathbb{D}}
\def\E{\mathcal{E}}
\def\O{\mathcal{O}}
\def\H{{\mathop{\rm H}}}

\def\P{\mathcal{P}}
\def\Q{\mathcal{Q}}

\def\M{\mathfrak{M}}
\def\N{\mathfrak{N}}

\def\p{{\sf{P}}}

\def\I{{\sf{I}}}
\def\J{{\sf{J}}}
\def\hotimes{\hat{\otimes}}
\def\hom{\mathcal{H}om}
\def\K{\mathbf{k}}
\def\proj{{\sf proj}}
\def\Proj{\mathbf{Proj}}

\def\mod{{\sf mod}}

\def\Vect{{\sf Vect}}
\def\Spec{\mathbf{Spec}}

\def\CM{\underline{\sf{CM}}}
\def\K{{\sf{K}}}
\def\Spc{\mathbf{Spc}}
\def\S{{\mathcal S}}
\def\D{{\sf D}}

\def\stmod{\underline{\sf mod}}

\def\R{{{\mathbb R}{\rm es}}}
\def\Inc{\mathbb{I}{\rm nc}}
\def\inc{{\mathop{\rm inc}}}
\def\supp{\mathbf{supp}}
\def\1{\mathbf{1}}

\title[]
 {Spectra of tensor triangulated categories over category algebras} 

\author{Fei Xu}

\address{Department of Mathematics, Shantou University, 243 University Road, Shantou, Guangdong 515063, China.}

\email{fxu@stu.edu.cn}

\thanks{The author \begin{CJK*}{UTF8}{}
\CJKtilde \CJKfamily{gbsn}徐 斐
\end{CJK*} was supported by a grant from the Department of Education of Guangdong Province.}


\begin{document}
\maketitle

\begin{abstract}
Let $\C$ be a finite EI category and $k$ be a field. We consider the category algebra $k\C$. Suppose $\K(\C)=\D^b(k\C\mbox{-}\mod)$ is the bounded derived category of finitely generated left modules. This is a tensor triangulated category and we compute its spectrum in the sense of Balmer. When $\C=G\propto\P$ is a finite transporter category, the category algebra becomes Gorenstein so we can define the stable module category $\CM k(G\propto\P)$, of maximal Cohen-Macaulay modules, as a quotient category of $\K(G\propto\P)$. Since $\CM k(G\propto\P)$ is also tensor triangulated, we compute its spectrum as well.

These spectra are used to classify tensor ideal thick subcategories of the corresponding tensor triangulated categories, despite the fact that the previously mentioned tensor categories are not rigid.
\end{abstract}

\section{Introduction}

Suppose $k$ is a field and $A$ is a $k$-algebra. We shall consider
left $A$-modules and denote by $A$-$\mod$ and $A$-$\proj$, respectively, the
categories of finitely generated $A$-modules and finitely generated
projective $A$-modules, respectively.

Suppose $\C$ is a finite category in the sense that $\Ob\C$ is a finite set. Consider the category algebra $k\C$. It is a finite-dimensional associative algebra with identity, and $k\C$-$\mod$ is a symmetric monoidal category, written as $(k\C\mbox{-}\mod,\hotimes,\k)$, in which the tensor identity $\k$ is also called the trivial $k\C$-module. Since the tensor product $-\hotimes-$ is exact in each variable, the bounded derived category of $k\C$-$\mod$, $\K(\C):=\D^b(k\C\mbox{-}\mod)$, is naturally a tensor triangulated category, although it is \textit{not} rigid in general.

In general, given a (rigid) tensor triangulated category $(\K,\otimes,\1)$, Balmer introduced an invariant of $\K$, called the triangular spectrum of $\K$, written as $\Spec\K=(\Spc\K,\O_{\Spc\K})$. Here $\Spc\K$, often referred to as the spectrum when there is no confusion, is the underlying topological space and $\O_{\Spc\K}$ is the structure sheaf of rings. The purpose of the present article is to compute $\Spc\K(\C)$. The structure sheaf is not easy to describe. However since $\O_{\Spc\K}$ is defined to be the sheafification of a presheaf $\O'_{\Spec\K}$, we can obtain some useful information about the latter.

When $\C$ is a finite EI category, $k\C$ enjoys some desirable homological properties. Particularly one can put a partial order on the set of isomorphism classes of objects in $\C$, which in turn allows a filtration of $k\C$-modules according to this partial order. Indeed a finitely generated $k\C$-module is built up by some modules of the automorphism groups in $\C$. When we compute the spectrum $\Spc\K(\C)$, the above-mentioned filtration provides a simple reduction to the computations of $\Spc\K(\Aut_{\C}(x))$. More precisely, there is a homeomorphism (Theorem 3.3.1)
$$
\Spc\K(\C)\simeq\biguplus_{[x]\subset\Ob\C}\Spc\K(\Aut_{\C}(x)),
$$
where $[x]$ runs over the set of isomorphism classes of objects in $\C$. Since $\Spc\K(G)$, for a finite group $G$, regarded as a category with one object, has been shown \cite{BCR2,BIK1,Ba1,Ba2} to be homeomorphic to $\Spec^h\H^{\cdot}(G,k)$, we obtain a satisfactory computational result on $\Spc\K(\C)$, for $\C$ finite EI. Here $\Spec^h\H^{\cdot}(G,k)$ stands for the spectrum of homogeneous prime ideals of $\H^*(G,k)$, the (finitely generated graded commutative) group cohomology ring.

We are particularly interested in finite transporter categories as they and their quotient categories are used widely in representation theory and algebraic topology. Let $G$ be a finite group and $\P$ be a finite $G$-poset. There exists a Grothendieck construction $G\propto\P$, which is a finite EI category and is commonly considered as the semi-direct product between $G$ and $\P$. Its category algebra $k(G\propto\P)$ is isomorphic to the skew group algebra $k\P[G]$, known to be Gorenstein. Recall that an algebra is called \textit{Gorenstein} if both the left and right regular modules are of finite injective dimension. Typical examples are self-injective algebras (e.g. group algebras $kG$) and algebras of finite global dimension (e.g. incidence algebras of finite posets $k\P$). For a Gorenstein algebra $A$, there exists a localizing sequence
$$
\D^b(A\mbox{-}\proj) \to \D^b(A\mbox{-}\mod) \to \CM(A),
$$
in which $\CM(A)\simeq \D^b(A\mbox{-}\mod)/\D^b(A\mbox{-}\proj)$ is
the \textit{stable module category} of $A$. When $A$ is self-injective,
$\CM(A)\simeq A$-$\stmod$ is the usual stable module category.
When $A$ is of finite global dimension, $\CM(A)\simeq\{0\}$. The category
$\CM(A)$ admits several equivalent characterizations, introduced
by Buchweitz \cite{B}. An equally illuminating (dual) construction can be found in Happel \cite{H}.

The category $\CM k(G\propto\P)$ is naturally tensor triangulated, so we can compute its spectrum $\Spc\CM k(G\propto\P)$. As a consequence of the homeomorphism
$$
\Spc\K(G\propto\P)\simeq\biguplus_{[x]\subset\Ob(G\propto\P)}\Spec^h\H^{\cdot}(\Aut_{G\propto\P}(x)),
$$
there exists another homeomorphism (Theorem 4.2.1)
$$
\Spc\CM k(G\propto\P)\simeq\biguplus_{[x]\subset\Ob(G\propto\P)}\Proj\H^{\cdot}(\Aut_{G\propto\P}(x)).
$$
The latter is used to classify tensor ideals of $\CM k(G\propto\P)$, while $\Spc\K(G\propto\P)$ can be used to classify tensor ideals of $\K(G\propto\P)$, by a theorem of Balmer (applied with care). These can be regarded as generalizations of certain results of Benson-Carlson-Rickard \cite{BCR2} and Benson-Iyengar-Krause \cite{BIK1}. The main feature is that we do not rely on a support variety theory. In fact, there exists an (ordinary) cohomology ring $\H^*(\C)$, which may be defined as the singular cohomology ring of the classifying space $B\C$. This graded commutative ring is frequently infinitely generated \cite{X3}, even modulo nilpotents. Thus we don't have a reasonable support variety theory, a widely adopted tool originated in group representations, which might have been used for understanding $k\C$-modules. In \cite{X1}, we identified finite transporter categories as a class of finite EI categories whose cohomology is finitely generated. Although there is a support variety theory \cite{X2}, we realized that the spectrum $\Spec^h\H^{\cdot}(G\propto\P)$ sometimes would not contain sufficient information to understand previously mentioned tensor triangulated categories (representation theory of $k(G\propto\P)$). For instance, when $\P$ is a connected poset (with trivial action by the trivial group), the cohomology ring $\H^*(\P)$ is finite-dimensional and thus $\Spec^h\H^{\cdot}(\P)$ is a point. By contrast, $k\P$-mod and $\D^b(k\P\mbox{-}\mod)$ have more than one tensor ideal thick subcategories, as long as $\P$ is not trivial. This observation led us to Balmer's theory, which turns out to solve our problem. In other words, we find in our cases that triangular spectrum seems to contain more representation-theoretic information than cohomology ring spectrum. We emphasize that these two kinds of spectra are closely related, as pointed out by Balmer. Indeed $\Spc\K(\C)$ may be considered as a topological space that affords all relevant cohomology rings in the form of sections of structure presheaf $\O'_{\Spc\K(\C)}$.

The article is organized as follows. We recall basics of transporter categories and category algebras in Section 2. Necessary ingredients of Balmer's tensor triangular geometry theory is also added for the convenience of the reader. Then in Section 3, we compute $\Spec\K(\C)$ of tensor triangulated categories, for a finite EI category $\C$. Various examples will be given to illustrate main results. Afterwards, Section 4 is devoted to study finite transporter categories.\\

\noindent {\bf Acknowledgements} I would like to thank Paul Balmer for explaining to me the possibility of using tensor triangular geometry theory to examine finite transporter categories. He read an earlier version of this article and kindly shared some ideas with me on the present work.

\section{Preliminaries}

We recall the definition of a category algebra and its basic properties. Particularly we shall describe the tensor structure on $k\C$-$\mod$ and $\K(\C)=\D^b(k\C\mbox{-}\mod)$. It is explained why $\K(\C)$ is not rigid. The definition of a transporter category is given afterwards. This section ends with a brief introduction to Balmer's tensor triangular geometry. 

Let us begin with the concept of a tensor triangulated category.

\begin{definition} A {\it tensor triangulated category} is a triple $(\K,\otimes,\1)$ consisting of a triangulated category $\K$, a symmetric monoidal (tensor) product $\otimes:\K\times\K\to\K$, which is exact in each variable and with respect to which there exists an identity $\1$.
\end{definition}

A {\it tensor triangular functor} $\mathbb{F} : \K \to \K'$ is an exact functor respecting the monoidal structures and preserves the tensor identity.

In a tensor triangulated category $\K$, we shall call a subcategory $\I$ a {\it tensor ideal} if it is a thick triangulated subcategory that is closed under tensoring with objects in $\K$.

\subsection{Tensor structure on $k\C$-mod and related tensor triangular categories}

Let us fix a field $k$ and a finite category $\C$ (that is, $\Mor\C$ is a finite set). Suppose $x\in\Ob\C$. We always write $[x]\subset\Ob\C$ for the set of objects that are isomorphic to $x$. 

The category algebra $k\C$ \cite{W, X} is defined as a $k$-vector space with basis the set of morphisms in $\C$. The multiplication is given on base elements by
$$
\beta * \alpha = \beta\circ\alpha
$$
if $\alpha : x \to y$ and $\beta : x' \to y'$ are composable (that is, $x'=y$), or zero otherwise. The category algebra $k\C$ is a finite-dimensional associative algebra with an identity
$1_{k\C}=\sum_{x\in\Ob\C}1_x$. Typical examples of finite categories are finite groups and finite posets. Their category algebras are the usual group algebras and incidence algebras. We emphasize that, however, there are many other interesting finite categories that are not of the two extreme cases, see for example \cite{X1} .

In this article we will often use the fact that a category equivalence $\mathcal{D}\to\C$ induces a Morita equivalence between their category
algebras, $k\mathcal{D}\simeq k\C$, as well as a homotopy equivalence between their classifying spaces $B\mathcal{D}\simeq B\C$ (see \cite{W}).

Let us denote by $\Vect_k$ the category of finite-dimensional $k$-vector spaces. It is a closed symmetric monoidal category equipped with a tensor product $-\otimes_k -$ and a tensor identity $k$. Suppose $\Vect_k^{\C}$ is the category of covariant functors from $\C$ to $\Vect_k$. Then a famous result of B. Mitchell says that
$$
\Vect^{\C}_k\simeq k\C\mbox{-}\mod.
$$
Since the symmetric monoidal structure on $\Vect_k$ naturally lifts to $\Vect^{\C}$, we obtain a symmetric monoidal category $(k\C\mbox{-}\mod,\hotimes_k,\k)$ (from now on we shall not distinguish $k\C$-$\mod$ from $\Vect^{\C}_k$). More precisely, the tensor product, $\hotimes_k$, is defined by
$$
(\M\hotimes_k \N) (x)=\M(x)\otimes_k \N(x)
$$
for any $\M, \N \in k\C$-$\mod$ and $x \in \Ob\C$. The module structure of $\M\hotimes\N$ can be viewed as given by the co-multiplication $\Delta : k\C \to k\C\otimes_k k\C$, induced by the canonical diagonal functor $\Delta : \C \to \C\times\C$ whose action on each $\alpha\in\Mor\C$ is $\Delta(\alpha)=\alpha\otimes\alpha$. The tensor product $-\hotimes_k-$ is exact in both variables. The constant functor $\k$, which takes the tensor identity $k$ of $\Vect_k$ as its value at each object, is the tensor identity with respect to $\hotimes_k$. Sometimes we also call $\k$ the trivial $k\C$-module because, when $\C$ is a group, $\k=k$ is exactly the trivial module of the group algebra. For the sake of simplicity, we shall write $\otimes$ for $\otimes_k$, and $\hotimes$ for $\hotimes_k$, throughout this article.

Since $-\hotimes-$ is exact in both variables, it gives rise to a tensor product on $\D^b(k\C\mbox{-}\mod)$. We shall still write $\hotimes$ and $\k$ for the tensor product and tensor identity in $\D^b(k\C\mbox{-}\mod)$.\\

\noindent {\bf Convention} We shall write $\K(\C)=\D^b(k\C\mbox{-}\mod)$. 

\subsection{EI categories}

When the category $\C$ is EI, in the sense that every endomorphism is an isomorphism, there is a good understanding of the representation theory of $k\C$ \cite{W}. Groups and posets are examples of EI categories. Given the EI property, the endomorphisms of each object form a finite group. Moreover, it allows us to introduce a natural partial order on the set of isomorphism classes of objects: $[x]\le[y]$ if and only if $\Hom_{\C}(x,y)\ne\emptyset$. This partial order in turn enables us to filtrate each $k\C$-module $\M$ by group modules. Indeed if $[x]$ is an isomorphism class of objects, then we may introduce a subspace $\M_x$ of $\M$ by $\M_x=\oplus_{y\in[x]}\M(y)$. For future reference, we write $\C_x$ for the full subcategory (a groupoid) of $\C$ consisting of objects $[x]$. Then $\M_x$ becomes a $k\C_x$-module. It is naturally a $k\C$-module too (not necessarily a submodule of $\M$).

When $[x]$ is maximal among the isomorphism classes of objects $[y]$ satisfying $\M(y)\ne 0$, $\M_x$ becomes a submodule of $\M$, and gives rise to a short exact sequence of $k\C$-modules
$$
0 \to \M_x \to \M \to \M/\M_x \to 0.
$$
Repeating the process we get a filtration on $\M$. When we examine $\K(\C)=\D^b(k\C\mbox{-}\mod)$, similarly each object $\c$ admits a filtration according to the same partial order.

Because of the filtration, every simple $k\C$-module $S$ must meet the condition that there exists one and only one isomorphism class of objects $[x]$ such that $S(x)\ne 0$. Moreover since $S(x)$ is a $k\Aut_{\C}(x)$-module, it has to be simple as a $k\Aut_{\C}(x)$-module. For future reference, the simple modules of $k\C$ are denoted by $S_{x,V}$, where $x$ stands for an isomorphism class $[x]\subset\Ob\C$ and $V=S_{x,V}(x)$ is a simple module of the automorphism group $\Aut_{\C}(x)$. Note that for each $[x]$, there is at least one simple module $S_{x,k}$. In general we have $\M_x\cong\M\hotimes S_{x,k}$. When $\C=\P$ is a poset, the automorphism groups are identities. It implies that $\{S_{x,k} \bigm{|} x\in\Ob\P\}$ is the set of all simple $k\P$-modules. Under the circumstance we shall write $S_x=S_{x,k}$ for brevity.

An element $\c\in\K(\C)$ is a bounded cochain complex of $k\C$-modules. Its cohomology $\H^*(\c)=\bigoplus_{i=m}^n\H^i(\c)$ is a direct sum of stalk complexes. Suppose $[x]\subset\Ob\C$. Then we may define an object of $\K(\C)$ as $\c_x\cong\c\hotimes S_{x,k}$. There is a K\"unneth formula \cite{X}, saying that $\H^*(\c\hotimes\d)=\H^*(\c)\hotimes\H^*(\d)$. Particularly as $k\C_x$-modules
$$
\H^*(\c_x)\cong\H^*(\c\hotimes S_{x,k})\cong\H^*(\c)\hotimes S_{x,k}.
$$

Recall in a tensor triangulated category $\K$ that an object $\c\in\K$ is called \textit{nilpotent} if $\c^{\otimes n}=0$ for some integer $n$. When $\C$ is EI, in $\K(\C)$ there is no such element. In fact if $\c^{\hotimes n}=0$, then $\H^*(\c^{\hotimes n})=\H^*(\c)^{\hotimes n}=0$. It forces $\H^i(\c)^{\hotimes n}=0$ for every $i$, and moreover
$$
\H^i(\c_x)^{\hotimes n}=[\H^i(\c)^{\hotimes n}]_x=0,
$$
for each $x\in\Ob\C$ (note that $S_{x,k}^{\hotimes n}=S_{x,k}$). Since $\H^i(\c_x)$ is a $k\C_x$-module, it forces $\H^i(\c_x)$=0 for every $x\in\Ob\C$. Hence $\H^*(\c)=0$, which implies $\c\cong 0$ in $\K(\C)$.

\begin{lemma} Let $\C$ be a finite EI category. Then there is no non-trivial nilpotent element in $\K(\C)$.
\end{lemma}

\subsection{$\K(\C)$ is not rigid}

Before we proceed, we want to point out that the category $\K(\C)$ is \textit{not} rigid. The reason to do so is that some statements in tensor triangular geometry \cite{Ba1,Ba2}, as well as in \cite{BCR2,BIK,BIK1}, requires the rigidity of a tensor triangular category. We need to stress on this subtlety, in order to avoid possible confusions.

By definition, a tensor triangular category $(\K,\otimes,\1)$ is \textit{rigid} if there exists a functor $D : \K\to\K^{op}$ satisfying
$$
\hom(\c,\d)\cong D(\c)\otimes\d
$$
for $\c,\d\in\K$. Here $\hom$ is the internal hom. If such a functor $D$ exists, then we must have $D(\c)\cong\hom(\c,\1)$. Thus when $D$ exists, we often call it a \textit{duality}. We shall focus on the case of $(\D^b(k\C\mbox{-mod}),\hotimes,\k)$, where $\C$ is a finite category.

First of all, let us recall the definition of the internal hom on $k\C$-mod and $\D^b(k\C\mbox{-}\mod)$ \cite{X}. Let $\M,\N \in k\C$-mod. Suppose $x\in\Ob\C$ and $k\C\cdot 1_x=k\Hom_{\C}(x,-)$ (a canonical projective $k\C$-module). We put
$$
\hom(\M,\N)(x):={\Hom}_{k\C}((k\C\cdot 1_x)\hotimes\M,\N).
$$
This functor $\hom(-,-)$ is left exact and hence induces an internal hom on the level of bounded derived category $\D^b(k\C\mbox{-}\mod)$. Let us still denote it by $\hom(-,-)$.

Suppose $\P$ is a finite poset. If $\M\in k\P$-$\mod$ and $x\in\Ob\P$, then we write $\M_{\ge x}$ for the submodule of $\M$ by ``brutal truncation''. By definition, we have
$$
\hom(\M,\N)(x)={\Hom}_{k\P}((k\P\cdot 1_x)\hotimes\M,\N)={\Hom}_{k\P}(\M_{\ge x},\N)\cong{\Hom}_{k\P}(\M_{\ge x},\N_{\ge x}).
$$

Next, we use a simple example to show that a duality functor $D : \K(\C) \to \K(\C)^{op}$ does not exist in general. Let $\P$ be the following poset
$$
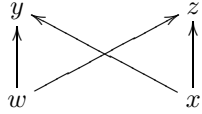
\xymatrix{y && z\\
w \ar[u] \ar[urr]&& x \ar[u] \ar[ull]}
$$
Denote by $S_x=S_{x,k}, S_y=S_{y,k}$ etc the 1-dimensional simple modules concentrated on designated vertices. Let $\M=S_y\oplus S_z$ be the radical of $\k$, that is, the kernel of $\k \twoheadrightarrow S_w\oplus S_x$, and $\N=S_w$. Then
$$
\hom(\M,\N)(w)={\Hom}_{k\P}(\M_{\ge w},S_{w})={\Hom}_{k\P}(S_y\oplus S_z,S_w)=0.
$$
If $D : k\P$-$\mod \to (k\P$-$\mod)^{op}$ were to exist, we would obtain
$$
[D(\M)\hotimes\N](w)=[\hom(\M,\k)\hotimes\N](w)={\Hom}_{k\P}(\M_{\ge w},\k)\otimes_k S_w(w)\cong k^2.
$$
By checking values on other objects, we find $\hom(\M,\N)=0$ while $D(\M)\hotimes\N\cong S_w^2$. This is a contradiction because there is no way that they are equal in $\K(\P)$. Thus, in general $\K(\P)=\D^b(k\P\mbox{-}\mod)$ is {\it not} rigid.

\subsection{Transporter categories}

A transporter category is a semi-direct product between a group and a poset. Directly from definition, it is EI. They appear in various places in representation theory and algebraic topology.

Let $G$ be a finite group, $\P$ a finite $G$-poset and $k$ a field (usually of positive characteristic $p\bigm{|}|G|$ in interesting cases). We will consider the transporter category $G\propto\P$, its category algebra $k(G\propto\P)$ and the bounded derived category
$$
\K(G\propto\P):=\D^b(k(G\propto\P)\mbox{-}\mod).
$$

We deem a group as a category with one object, usually denoted by
$\bullet$. The identity of a group $G$ is always named $e$. We say a
poset $\P$ is a $G$-poset if there exists a functor $F$ from $G$ to
$\mathfrak{sCats}$, the category of small categories, such that
$F(\bullet)=\P$. It is equivalent to saying that we have a group
homomorphism $G \to \Aut(\P)$.

\begin{definition} Let $G$ be a group and $\P$ a $G$-poset. The
transporter category $G\propto \P$ has the same objects as $\P$,
that is, $\Ob(G\propto \P)=\Ob\P$. For $x ,y \in \Ob(G\propto \P)$,
the morphism set $\Hom_{G\propto\P}(x,y)$ is
$$
\{(g, gx\le y) \bigm{|} g\in G\}.
$$

If $(g, gx\le y)$ and $(h,hy\le z)$ are two morphisms in
$G\propto\P$, then their composite is easily seen to be $(h,hy\le z)\circ(g, gx\le y)=(hg,
(hg)x\le z)$.
\end{definition}

From the definition one can easily see that there is a natural
embedding $\iota : \P\hookrightarrow G\propto\P$ via $(x \le y)
\mapsto (e, x\le y)$. On the other hand, the transporter category
admits another natural functor $\pi : G\propto\P \to G$, given by
$x\mapsto \bullet$ and $(g,gx\le y)\mapsto g$. Thus we always have a
sequence of functors
$$
\P {\buildrel{\iota}\over{\hookrightarrow}} G\propto\P
{\buildrel{\pi}\over{\longrightarrow}} G
$$
such that $\pi\circ\iota(\P)$ is the trivial subgroup or subcategory
of $G$. Topologically it is well known that $B(G\propto\P)\simeq EG\times_G B\P$. Passing to classifying
spaces, we obtain a fibration sequence
$$
B\P {\buildrel{B\iota}\over{\longrightarrow}} EG\times_G B\P
{\buildrel{B\pi}\over{\longrightarrow}} BG .
$$
Forming the transporter category over a $G$-poset eliminates the
$G$-action, and thus is the algebraic analogy of introducing the
Borel construction over a $G$-space. As a consequence, we have
$$
\H^*(G\propto\P):=\Ext^*_{k(G\propto\P)}(\k,\k)\cong\H^*(EG\times_GB\P,k)=\H^*_G(B\P,k).
$$
Here $\Ext^*_{k(G\propto\P)}(\k,\k)$ is called the \textit{ordinary cohomology ring} of $k(G\propto\P)$, and $\H^*_G(B\P,k)$ is the \textit{equivariant cohomology ring} of the $G$-space $B\P$. This is a finitely generated graded commutative ring \cite{X1}.

We note that for each $x\in\Ob(G\propto\P)=\Ob\P$, $\Aut_{G\propto\P}(x)$ is exactly
the isotropy group of $x$. For the sake of simplicity, we will often
write $G_x=\Aut_{G\propto\P}(x)$, and $[x]=G.\{x\}$, the orbit of $x$. Note that
$[x]$ is a $G$-(po)set, consisting of exactly the objects in $G\propto\P$ that
are isomorphic to $x$.

\begin{example}
\begin{enumerate}
\item If $G$ acts trivially on $\P$, then $G\propto\P = G\times \P$. In this case for any $x\in\Ob(G\times\P)$, $G_x=G$. Here $\times$ stands for the product between two categories.

\item Let $G$ be a finite group and $H$ a subgroup. We consider
the set of left cosets $G/H$ which can be regarded as a $G$-poset:
$G$ acts via left multiplication. The transporter category $G\propto
(G/H)$ is a connected groupoid whose skeleton is isomorphic to $H$.
It means there is no essential
difference between $H$ and $G\propto(G/H)$ as far as we concern.
Hence it makes sense if we regard transporter categories as
generalized subgroups for a fixed finite group.

For an arbitrary $G$-poset $\P$ and $x\in\Ob\P$, we have a category
equivalence $G\propto[x]=G\propto G.\{x\}\simeq G_x$.
\end{enumerate}
\end{example}

\noindent {\bf Convention} When $G=\{e\}$ is trivial, we simplify $\K(\{e\}\propto\P)=\K(\{e\}\times\P)$ to $\K(\P)$. While when $\P=\bullet$ is trivial, we denote by $\K(G)$ the category $\K(G\propto\bullet)=\K(G\times\bullet)=\D^b(kG\mbox{-}\mod)$.

\subsection{Tensor triangular geometry}

Let $\K$ be an essentially small tensor triangulated category, written as $(\K, \otimes, \1)$. We record basics of Balmer's tensor triangular geometry theory \cite{Ba1,Ba2}.

\begin{definition} Suppose $\p$ is a tensor ideal of $\K$. It is said to be \textit{prime} if $\p$ is properly contained in $\K$ and if $\c\otimes\d\in\p$ implies either $\c\in\p$ or $\d\in\p$. Denote by $\Spc\K$ the set of all prime ideals of $\K$. If $\c\in\K$, its \textit{support} is defined to be
$$
\supp(\c)=\{\p \in \Spc\K \bigm{|} \c\not\in\p\}.
$$
One can topologize $\Spc\K$ by asking the following to be an open basis
$$
U(\c)=\Spc\K-\supp(\c)=\{\p \in \Spc\K \bigm{|} \c\in\p\}.
$$

Indeed every quasi-compact open subset of $\Spc\K$ is of the form $U(\c)$ for some $\c\in\K$. One can continue to introduce a \textit{structure sheaf} of rings on $\Spc\K$. Let $U$ be a quasi-compact open subset of $\Spc\K$ and $Z=\Spc\K-U$ its complement. Define a tensor ideal of $\K$ as
$$
\K_Z=\{\c\in\K\bigm{|}\supp(\c)\subset Z\}.
$$
Then one gets a tensor triangulated category
$$
\K(U)=(\K/\K_Z)^{\natural},
$$
as the idempotent completion of $\K/\K_Z$. A presheaf $\O'_{\Spc\K}$ of rings on $\Spc\K$ is given by
$$
\O'_{\Spc\K}(U)=\End_{\K(U)}(\1),
$$
and $\O_{\Spc\K}$ is defined to be its sheafification. Here $\1$ is the tensor identity of $\K(U)$, as the canonical image of $\1\in\K$. In this way $\Spec\K=(\Spc\K,\O_{\Spc\K})$ becomes a locally ringed space, called the \textit{spectrum} of $\K$. Sometimes one abuses the terminology and calls the underlying space $\Spc\K$ the spectrum of $\K$.
\end{definition}

It was shown by Balmer that $\Spc\K$ is not empty as long as $\K\ne \{0\}$.

\begin{definition} If $\I\subset\K$   is a tensor ideal in a tensor triangulated category $(\K,\otimes,\1)$, then its radical is defined to be
$$
\sqrt{\I}=\{\c\in\K\bigm{|} \exists  n\ge 1\ \mbox{such that}\ \c^{\otimes n}\in\I\}.
$$
If $\I=\sqrt{\I}$, we say this tensor ideal is radical.
\end{definition}

We will refer to the following result a couple of times.

\begin{theorem}[Balmer] Suppose $\Spc\K$ is a Noetherian topological space. Then there is a bijection between the set of specialization closed subsets of $\Spc\K$ and the set of radical tensor ideals of $\K$.
\end{theorem}

\begin{remark} It is very often in a tensor triangulated category $\K$ that all tensor ideals are radical. For instance $\K(G)$, for a finite group $G$, is one of the examples. From here we can deduce in Section 3.1 that all tensor ideals of $\K(\C)$ are radical.
\end{remark}

Every tensor triangulated functor $F : \K\to\K'$ induces a continuous map $\Spc(F):\Spc\K'\to\Spc\K$. Some other constructions and results of Balmer will be quoted later on when they are needed.

\section{Spectrum of $\K(\C)$}

From now on, we always assume $\C$ to be a finite EI category. We shall apply tensor triangular geometry theory to investigate $\K(\C)=\D^b(k\C\mbox{-}\mod)$. One may bear the simplest case, $\K(\P)$, of a poset in mind as a prototype.

There exists an analogy for quiver representations \cite{LS}. 

\subsection{The restriction}

Suppose $\E$ is a full subcategory of $\C$. Then the inclusion $\E \hookrightarrow \C$ induces a restriction (brutal truncation)
$$
\res_{\E}=\res^{\C}_{\E} : k\C\mbox{-}\mod \to k\E\mbox{-}\mod.
$$
It is exact and preserves both tensor products and tensor identity. We write the resulting tensor derived functor as
$$
\R_{\E}=\R^{\C}_{\E} : \K(\C) \to \K(\E),
$$
which gives rise to a continuous map
$$
\Spc(\R_{\E}) : \Spc\K(\E) \to \Spc\K(\C).
$$

\begin{definition}
Let $\C$ be an EI category. A full subcategory $\E$ is said to be {\it convex} if, given objects $x, z \in \Ob\E$ and $y\in\Ob\C$, we must have $y\in\Ob\E$ whenever there exist morphisms $x \to y$ and $y\to z$.
\end{definition}

If $\E\subset\C$ is convex, then every $k\E$-module is automatically a $k\C$-module. Thus we have a fully faithful functor
$$
\inc_{\E}=\inc^{\C}_{\E}: k\E\mbox{-}\mod \to k\C\mbox{-}\mod
$$
It is also exact, preserving tensor products and tensor identity. We write the resulting tensor derived functor as
$$
\Inc_{\E}=\Inc^{\C}_{\E} : \K(\E) \to \K(\C),
$$
which also gives rise to a continuous map
$$
\Spc(\Inc_{\E}) : \Spc\K(\C) \to \Spc\K(\E).
$$

Let $x\in\Ob\C$. Then $\C_{x}$ is a convex subcategory. We note that for a tensor ideal $\I\subset\K(\C)$, $\R_{\C_x}\I$ is a tensor ideal of $\K(\C_x)$. Since
$$
\R_{\C_x}\I = S_{x,k}\hotimes\I,
$$
$S_{x,k}\hotimes\I=\Inc_{\C_x}(S_{x,k}\hotimes\I)\subset\I$ 
can also be regarded as a tensor ideal of $\K(\C)$.

\begin{lemma}  Suppose $\E\subset\C$ is a full subcategory. We write $\C\backslash\E$ for the full subcategory consisting of objects $x\not\in\Ob\E$.
\begin{enumerate}
\item Suppose $\E\subset\C$ is a convex subcategory. Then $\R_{\E}\circ\Inc_{\E}=\Id$ and $\R_{\E}\circ\Inc_{\E}=0$.

\item Suppose $\E\subset\C$ is an arbitrary full subcategory. Then $\ker(\R_{\C\backslash\E})$ is generated over $\Inc_{\C_x}\K(\C_x)$, where $x$ runs over $\Ob\E$.
\end{enumerate}

\begin{proof} For (1), we observe that $\res_{\E}\circ\inc_{\E}=\Id$ and $\res_{\C\backslash\E}\circ\inc_{\E}=0$. The equalities follow from these. The first equation implies that $\inc_{\E}$, hence $\Inc_{\E}$, are fully faithful.

As for (2), given $x\in\Ob\C$, $\C_x$ is a convex subcategory. Thus $\Inc_{\C_x}\K(\C_x)$ becomes a tensor subcategory of $\K(\C)$. Since $x\in\Ob\E$ and $\C\backslash\E \subset \C\backslash\C_x$, we have $\R_{\C\backslash\E}\circ\Inc_{\C_x}\K(\C_x)=\R_{\C\backslash\E}\circ\Inc_{\E}\K(\C_x)=0$, which implies $\Inc_{\C_x}\K(\C_x) \subset \ker(\R_{\C\backslash\E})$. On the other hand, if we take $y\not\in\Ob\E$, then $\R_{\C\backslash\E}\circ\Inc_{\C_y}(M)=M$ for every $k\C_y$-module $M$ regarded as a $k\C$-module. It means $\Inc_{\C_y}\K(\C_y)\cap\R_{\C\backslash\E}=\{0\}$ because $\Inc_{\C_y}$ is fully faithful. Since $\Inc_{\C_x}\K(\C_x)$, $x\in\Ob\C$, generate $\K(\C)$, we have proved the claim.
\end{proof}
\end{lemma}

The kernel of a tensor triangulated functor is always a tensor ideal.

\begin{example} Let $\P$ be a poset and $\Q$ a convex subposet. Then $\res_{\Q}: k\P$-$\mod \to k\Q$-$\mod$ and $\inc_{\Q} : k\Q$-$\mod \to k\P$-$\mod$ are easy to write down. Consequently $\R_{\Q}$ sends $\K(\P)$ onto $\K(\Q)$ while $\Inc_{\Q}$ maps $\K(\Q)$ into $\K(\P)$ with image equal to
$$
\langle S_x=S_{x,k} \bigm{|} x \in \Ob\Q \rangle,
$$
the tensor ideal generated by $S_x, x\in\Ob\Q$.
\end{example}

\begin{lemma} The functor $\R_{\E}$ preserves tensor ideals and if $\E$ is convex, so does $\Inc_{\E}$. Hence if $\E$ is convex and $\I\subset\K(\C)$ is a tensor ideal, so is $\Inc_{\E}\circ\R_{\E}(\I)$.

\begin{proof} For $\R_{\E}$, it follows from the fact that it is full and surjective on objects. 

Assume $\E$ is convex. Suppose $\J\subset\K(\E)$ is a tensor ideal. Let $\c\in\K(\C)$ and $\d\in\J$. Then $\c\hotimes\Inc_{\E}(\d)\cong[\Inc_{\E}\circ\R_{\E}(\c)]\hotimes\Inc_{\E}(\d)\cong\Inc_{\E}(\R_{\E}(\c)\hotimes\d)\in\Inc_{\E}(\J)$. It implies that $\Inc_{\E}(\J)$ is a tensor ideal. Particularly $\Inc_{\E}\circ\R_{\E}(\I)$ is a tensor ideal if $\I\subset\K(\C)$ is.

The containment $\Inc_{\E}\circ\R_{\E}(\I)\subset\I$ follows from the fact that if $\c\in\I$, then $\c_x=\c\hotimes S_{x,k}\in\I$ for every $x\in\Ob\C$. It implies that $\Inc_{\E}\circ\R_{\E}(\c)\in\I$ because $\c$ admits a filtration by these $\c_x$.
\end{proof}
\end{lemma}

\begin{example} Let $\P$ be a poset. Suppose $\I\subset\K(\P)$ is a tensor ideal. We define a tensor ideal of $\K(\P)$ as
$$
\I'=\langle S_x\bigm{|} S_x\hotimes\c\neq 0\ \mbox{for some}\ \c\in\I\rangle.
$$
Since the defining condition $0 \neq S_x\hotimes\c=\c_x\cong\H^*(\c_x)$ (by the splitting of complexes of $k$-vector spaces) implies $S_x\in \I$, the above tensor ideal is a tensor subcategory of $\I$. It follows immediately that $\I'=\I$. Furthermore we have $\I=\ker(\R_{\P\backslash\Q})$, where $\Q=\{x\bigm{|}S_x\in\I\}$ is a subposet. Thus every tensor ideal of $\K(\P)$ is of the form $\ker(\R_{\Q'})$ for some subposet $\Q'\subset\P$.
\end{example}

\begin{example} Let $\P = x \to y \to z$. Each object determines a tensor ideal so we can consider $\ker(\R_{y\to z}),\ker(\R_{x\to z}),\ker(\R_{x\to y})$. Note that the intersection between any two of them becomes zero. These tensor ideals are not prime. If we take $S_x$ and $S_y$ to be the 1-dimensional simple modules supported on $x$ and $y$, respectively, then regarding them as stalk complexes we have $S_x \in \ker(\R_{y\to z})$ and $S_y\in\ker(\R_{x\to z})$. But $S_x\hotimes S_y = 0 \in\ker(\R_{x\to y})$ excludes $\ker(\R_{x\to y})$ from being prime. The same argument tells that $\ker(\R_{y\to z}), \ker(\R_{x\to z}), \{0\}=\ker(\Id)$ are not prime.

We will see in short that the rest proper tensor ideals $\ker(\R_{\P_z})$, $\ker(\R_{\P_x})$ and $\ker(\R_{\P_y})$ ($\P\backslash\P_y=x\to z$ is not convex!) are the only prime ones in $\K(\P)$. We will also see that surprisingly the tensor ideals are determined only by $\Ob\P$. It means we can forget about the non-isomorphisms (visible arrows).
\end{example}

\begin{lemma} Suppose $\I\subset\K(\C)$ is a tensor ideal. Then it is radical.

\begin{proof} In fact if $\c\in\K(\C)$ and $\c^{\hotimes n}\in\I$ for some integer $n$, then for every $[x]\subset\Ob\C$ we have $\R_{\C_x}(\c^{\hotimes n})\in\R_{\C_x}\I$. Since $\R_{\C_x}(\c^{\hotimes n})=(\Res_{\C_x}\c)^{\hotimes n}$ and $\R_{\C_x}\I\subset\R_{\C_x}\K(\C)=\K(\C_x)\simeq\K(\Aut_{\C}(x))$ is a tensor ideal, by Remark 2.5.4 we know $\c_x=\Res_{\C_x}\c\in\R_{\C_x}\I=S_{x,k}\hotimes\I\subset\I$ for every $[x]\subset\Ob\C$. It follows that $\c\in\I$.
\end{proof}
\end{lemma}

\subsection{Prime ideals}

Consider the functor $\R_{\E} : \K(\C) \to \K(\E)$. If $\I$ is a tensor (resp. prime) ideal of $\K(\E)$, then $\R_{\E}^{-1}(\I)$ is a tensor (resp. prime) ideal of $\K(\C)$ \cite{Ba1,Ba2}.

Recall that we denote by $[x]$ the isomorphism class of $x$ in $\C$. Filtration on complexes of $k\C$-modules gives us a characterization of prime ideals.

\begin{lemma} Given a convex subcategory $\E\subset\C$ and a tensor ideal $\I\subset\K(\E)$, the preimage $\R_{\E}^{-1}(\I)$ can be identified with
$$
\langle \Inc_{\E}(\I), \ker(\R_{\E}) \rangle=\langle \Inc_{\E}(\I), S_{x,V} \bigm{|} x\not\in\Ob\E \rangle.
$$

\begin{proof} Consider the composite of functors $\K(\C) {\buildrel{\R_{\E}}\over{\to}} \K(\E) \to \K(\E)/\I$. Use the fact that $\R_{\E}\circ\Inc_{\E}=\Id$.
\end{proof}
\end{lemma}

Particularly if $\p\in\K(\C_x)$ is a prime ideal, then its preimage is the prime ideal
$$
\R_{\C_x}^{-1}(\p)=\langle \Inc_{\C_x}(\p),\ker(\R_{\C_x})\rangle.
$$

\begin{lemma} Let $\p\subset\K(\C)$ be a prime ideal. Then there is exactly one $[x]\subset\Ob\C$ such that $\R_{\C_x}\p\varsubsetneq\R_{\C_x}\K(\C)$, whence $\p=\langle \Inc_{\C_x}\circ\R_{\C_x}\p,\ker(\R_{\C_x})\rangle$.

\begin{proof} If there are $x\not\cong y$ such that both $\R_{\C_x}\p\varsubsetneq\R_{\C_x}\K(\C)$ and $\R_{\C_y}\p\varsubsetneq\R_{\C_y}\K(\C)$. Then there exist $S_{x,V}\not\in\Inc_{\C_x}\circ\R_{\C_x}\p$ and $S_{y,W}\not\in\Inc_{\C_y}\circ\R_{\C_y}\p$, because the simple modules of $k\C_z$ generate $\K(\C_z)=\R_{\C_z}\K(\C)$ for each isomorphism class $[z]\subset\Ob\C$. Since neither $S_{x,V}$ nor $S_{y,W}$ lies in $\p$, the fact that $S_{x,V}\hotimes S_{y,W}=0\in\p$ provides a contradiction to the assumption that $\p$ is prime.

There has to be such an isomorphism class $[x]$ satisfying the condition that $\R_{\C_x}\p\varsubsetneq\R_{\C_x}\K(\C)$, because otherwise $\p$ is not properly contained in $\K(\C)$.
\end{proof}
\end{lemma}

The next result continue to characterize $\R_{\C_x}\P$, when $\P$ is prime.

\begin{lemma} Let $\p\in\K(\C)$ be a prime ideal and $\E\subset\C$ a convex subcategory. Then $\R_{\E}(\p)$ is either $\K(\E)$ or a prime ideal in $\K(\E)$. 

\begin{proof} Assume $\R_{\E}\P\neq\K(\E)$. Suppose $\c,\d\in\K(\E)$ and $\c\hotimes\d\in\R_{\E}(\p)$. Since $\Inc_{\E}(\c\hotimes\d)=\Inc_{\E}(\c)\hotimes\Inc_{\E}(\d)\in\Inc_{\E}\circ\R_{\E}\p\subset\p$, it follows that either $\Inc_{\E}(\c)$ or $\Inc_{\E}(\d)$ belongs to $\p$. Thus either $\c=\R_{\E}\circ\Inc_{\E}(\c)$ or $\d=\R_{\E}\circ\Inc_{\E}(\d)$ belongs to $\R_{\E}(\p)$. It means that $\R_{\E}\P$ is prime.
\end{proof}
\end{lemma}

In summary we have the following characterization of prime ideals. Note that when $\C$ has more than one isomorphism classes of objects then $\{0\}$ is not a prime ideal (think about the simplest case of a poset).

\begin{proposition} Let $\p\subset\K(\C)$ be a prime ideal. Then  $\R_{\C_x}(\p)\neq\K(\C_x)$ for a unique $[x]$. Whence $\R_{\C_x}(\p)$ is a prime ideal of $\K(\C_x)$), and
$$
\p=\R_{\C_x}^{-1}(\R_{\C_x}(\p))=\langle \Inc_{\C_x}\circ\R_{\C_x}\p,\ker(\R_{\C_x})\rangle.
$$
\end{proposition}

\begin{example} Let $\P$ be a poset and $x\in\Ob\P$. The object determines a subposet $\P_x$. The restriction $\res_{\P_x}$, also called the evaluation at $x$, $\res_{\P_x} : k\P$-$\mod \to k\P_x$-$\mod\cong\Vect_k$, gives rise to a tensor triangulated functor $\R_{\P_x} : \K(\P) \to \K(\P_x)=\D^b(\Vect_k)$, whose kernel is $\langle S_y \bigm{|} y\ne x\rangle$. If we write $\P\backslash\{x\}$ as the maximal subposet of $\P$ with $x$ removed, then the kernel $\ker(\R_{\P_x})$ is maximal and hence is prime in $\K(\P)$. When $x$ runs over the set of (isomorphism classes of) objects, they give us all possible maximal ideals in $\K(\P)$.

Let $\I$ be a prime ideal in $\K(\P)$. Then $\I\subset\p^x$ for some $x\in\Ob\P$. Since $S_y\hotimes S_x=0\in\I$ and $S_x\not\in\I$, by the definition of a prime ideal, $S_y\in\I$ for all $y\ne x$. It means that $\I=\p^x$.

The spectrum $\Spc\K(\P)=\{\p^x\bigm{|}x\in\Ob\P\}$ has finitely many points and thus has a simple topological structure. Given $x\in\Ob\P$,
$$
\supp(S_x)=Z(\{S_x\})=\{\p\in\Spc\K(\P) \bigm{|}\{S_x\}\cap \p=\emptyset\}=\{\p^x\}
$$
is a closed subset. It means $\Spc\K(\P)$ is discrete. More generally if $\c\in\K(\P)$, then
$$
\supp(\c)=\{\p\in\Spc\K(\P) \bigm{|}\{\c\}\cap \p=\emptyset\}=\{\p^x\bigm{|}S_x\hotimes\c\ne 0, x\in\Ob\P\}.
$$
\end{example}

\subsection{Triangular spectrum}

Back to the general situation of finite EI categories, we prove that the spectrum $\Spc\K(\C)$ is the disjoint union of $\Spc\K(\C_x)\simeq\Spc\K(\Aut_{\C}(x))$, $[x]$ running over the set of isomorphism classes of objects in $\C$.

\begin{theorem} There are homeomorphisms
$$
\Spc\K(\C)\simeq\biguplus_{[x]\subset\Ob\C}\Spc\K(\C_x)\simeq\biguplus_{[x]\subset\Ob\C}\Spec^h\H^{\cdot}(\Aut_{\C}(x)),
$$
where $\Spec^h\H^{\cdot}(\Aut_{\C}(x))$ is the spectrum of homogeneous prime ideals of the group cohomology ring $\H^{\cdot}(\Aut_{\C}(x))$.

\begin{proof} The second homeomorphism is true because $\C_x\simeq\Aut_{\C}(x)$ as categories. It means that $k\C_x$ is Morita equivalent to $k\Aut_{\C}(x)$. Then by \cite{Ba2,BIK1}, there is a homeomorphism $\Spc\K(\Aut_{\C}(x))\simeq\Spec^h\H^{\cdot}(\Aut_{\C}(x))$. We continue to prove the first homeomorphism. In fact there is an obvious bijection that provides a set isomorphism between them. Suppose $\p\in\Spc\K(\C)$. We define a map
$$
\widetilde{\R} : \Spc\K(\C) \to \biguplus_{[x]\subset\Ob\C}\Spc\K(\C_x)
$$
by $\p \mapsto \R_{\C_x}(\p)$, where $[x]$ is the only $G$-orbit such that $\R_{\C_x}(\p)\varsubsetneq\R_{\C_x}(\K(\C))$. We shall prove this bijection and its inverse are continuous.

Suppose $\c\in\K(\C)$. We can compute its support as
$$
\begin{array}{ll}
\supp(\c)&=\{\p\in\Spc\K(\C)\bigm{|}\c\not\in\p\}\\
&\\
&=\{\R_{\C_x}^{-1}(\R_{\C_x}(\p))\in\Spc\K(\C)\bigm{|}\R_{\C_x}(\c)\not\in\R_{\C_x}(\p), [x]\subset\Ob\C\}\\
&\\
&=\biguplus_{[x]\subset\Ob\C}\{\R_{\C_x}^{-1}(\R_{\C_x}(\p))\in\Spc\K(\C)\bigm{|}\R_{\C_x}(\c)\not\in\R_{\C_x}(\p)\}\\
&\\
&\leftrightarrow\biguplus_{[x]\subset\Ob\C}\{\R_{\C_x}(\p)\in\Spc\K(\C_x)\bigm{|} \R_{\C_x}(\c)\not\in\R_{\C_x}(\p)\}.
\end{array}
$$
Thus the map $\widetilde{\R}$ restricts to bijections on closed subsets of the two spaces $\Spc\K(\C)$ and $\bigcup_{[x]\subset\Ob\C}\Spc\K(\C_x)$:
$$
\supp(\c)\mapsto\biguplus_{[x]\subset\Ob\C}\supp(\R_{\C_x}(\c)).
$$
We are done.
\end{proof}
\end{theorem}

As an example when $G$ is a finite group, $\P=\S_p$ is the poset of non-trivial $p$-subgroups of $G$ and $\C=\O_p(G)$ is the $p$-orbit category, the automorphism groups of $\O_p(G)$ are $N_G(P)/P$, for every $P\in\Ob\O_p(G)$. Thus
$$
\Spc\K(\O_p(G))\simeq\biguplus_{[P]\subset\Ob\O_p(G)}\Spec^h\H^{\cdot}(N_G(P)/P).
$$

The structure sheaf on $\Spc\K(\C)$ is usually hard to compute. But various relevant cohomology rings appear as the sections of the presheaf $\O'_{\Spc\K(\C)}$. Indeed to each simple $k\C$-module $S_{x,k}$, we have
$$
\supp(S_{x,k})=\{\p\in\Spc\K(\C)\bigm{|}\Res_{\C_x}\p=\K(\C_x)\}.
$$
Since $\supp(S_{x,k}\oplus S_{y,k})=\supp(S_{x,k})\cup\supp(S_{y,k})$, given a full subcategory $\E$ we can take an open subset $U_{\E}=\Spc\K(\C)-Z_{\E}$ of $\Spc\K(\C)$, where $Z_{\E}=\biguplus_{[x]\not\subset\Ob\E}\supp(S_{x,k})$, a closed subset of $\Spc\K(\C)$. Then $U_{\E}=\bigcap_{[x]\not\subset\Ob\E}U(S_{x,k})$ is quasi-compact. It turns out that $\K(\C)(U_{\E})=(\K(\C)/\K(\C)_{Z_{\E}})^{\natural}\simeq\K(\E)$, and that
$$
\O'_{\Spc\K(\C)}(U_{\E})=\End_{\K(\C)(U_{\E})}(\k)\cong\H^*(\E).
$$
To some extent, the underlying space $\Spc\K(\C)$ only sees the isomorphism classes of objects of $\Ob\C$, while along with the presheaf, it reveals deep structural information.

\begin{lemma} Let $\P$ be a finite poset. Then the sections of $\O'_{\Spc\K(\P)}$ are exactly the cohomology rings of all subposets of $\P$. The structure sheaf $\O_{\Spc\K(\P)}$, a.k.a. the shifification of $\O'_{\Spc\K(\P)}$, is constant.

\begin{proof} Let $U\subset\Spc\K(\P)$ be an open set. It is always quasi-compact since $\Spc\K(\P)$ is a finite space. Suppose $Z=\{\p^x=\ker(\R_{\P_x})\bigm{|} x\in\Ob\Q\}$, where $\Q$ is a subposet of $\P$, is the complement of $U$. Then
$$
\K(\P)_Z=\{a\in\K(\P)\bigm{|}\supp(a)\subset Z\}=\langle S_x \bigm{|} x\in\Ob\Q\rangle=\ker(\R_{\P\backslash\Q}).
$$
It means $\K(\P)(U)=[\K(\P)/\K(\P)_Z]^{\natural}\simeq\K(\P\backslash\Q)$. Thus $\O'_{\Spc\K(\P)}(U)\cong\End_{\K(\P\backslash\Q)}(\k)\cong\H^*(\P\backslash\Q)$ is the cohomology ring of the poset $\P\backslash\Q$ with coefficients in $k$.

To understand the sheafification of $\O'_{\Spc\K(\P)}$, we need to compute the stalks at points of $\Spc\K(\P)$. Because $\Spc\K(\P)$ is discrete, each point $\p^x$ forms an open neighborhood of itself. Hence the stalk at $\p^x$ is $[\O'_{\Spc\K(\P)}]_{\p^x}=[\O'_{\Spc\K(\P)}](\p^x)\cong\H^*(\{x\})=k$, and this implies that $\O_{\Spc\K(\P)}$ is constant.
\end{proof}
\end{lemma}

\begin{remark}
The space $\Spc\K(-)$ alone, as an invariant, is not stronger than $\H^*(-)$, because for example it does not distinguish the two-object posets $x\ \cdot\to\cdot\ y$ and $x\cdot \ \ \ \ \cdot\ y$. This is also reflected in the fact that Balmer's canonical map
$$
\rho^{\cdot} : \Spc\K(\P) \to \Spec^h(\H^{\cdot}(\P))
$$
is surjective, but not homeomorphic in general. However if we take the presheaf of rings $\O'_{\Spc\K(\P)}$ on $\Spc\K(\P)$, and note that its global section is $\End_{\K(\P)}(\k)\cong\H^*(\P)$, then we see $\O'_{\Spc\K(-)}$ does distinguish $\P$ and $\Q$.

In the general situation, the structure presheaf on $\Spc\K(\C)$ is determined to a very large extent by the intrinsic structure of $\C$. It makes $\Spec\K(\C)$ a more powerful invariant than $\biguplus_{[x]\subset\Ob\C}\Spec^h\H^{\cdot}(\Aut_{\C}(x))$.
\end{remark}

Since $\Spec^h\H^{\cdot}(G)$ is Noetherian for any finite group $G$, so is $\Spc\K(\C)$ as a union of finitely many Noetherian spaces.

\begin{corollary} There is a bijection between the set of tensor ideals in $\Spc\K(\C)$ and the set of specialization closed subsets of $\Spc\K(\C)$.
\end{corollary}

Let us work out the classification for posets.

\begin{example} Suppose $\I\subset\K(\P)$ is a tensor ideal. Assume $S_x\hotimes\I\ne 0$. Then, because of the splitting of exact sequences of $k$-vector spaces, we deduce that $S_x\in\I$.  Thus $\I=\langle S_x\bigm{|} S_x\hotimes\I\ne 0 \rangle$. In other words, $\I$ is determined by a set of simple $k\P$-modules, or equivalently, a subset of $\Ob\P$ or a subposet of $\P$. Suppose $\Q_{\I}\subset\P$ corresponds to $\I$. The poset $\Q_{\I}$ is the intersection of a certain set of maximal subposets
$$
\{\P\backslash\{x\}\bigm{|}x\not\in\Ob\Q_{\I}\}.
$$
We aim to establish a one-to-one correspondence between the set of subsets of $\Spc\K(\P)$ and the set of tensor ideals of $\K(\P)$.

Indeed we have the following bijections
$$
\xymatrix{\Ob\P \ar@{<->}[r] & \{\mbox{maximal subposets of}\ \P\} \ar@{<->}[r] & \{\mbox{prime ideals of}\ \K(\P)\}=\Spc\K(\P),\\
x \ar@{|-{>}}[r] & \P\backslash\{x\} \ar@{|-{>}}[r] & \p^x=\ker(\R_{\P_x})}
$$
which can be naturally extended to some 1-to-1 correspondences
$$
\xymatrix{\{\mbox{subsets of}\ \Ob\P\} \ar@{<->}[rr] \ar@{<->}[d] && \{\mbox{subsets of}\ \Spc\K(\P)\} \ar@{<->}[d]\\
\{\mbox{subposets of}\ \P\} \ar@{<->}[rr] && \{\mbox{tensor ideals of}\ \K(\P)\}}
$$
In the preceding diagram, the definitions of the left vertical correspondence and the lower horizontal one are known. To describe the top horizontal correspondence, let $Z\subset\Spc\K(\P)$. Then $Z=\{\p^x\bigm{|}\mbox{some}\ x\in\Ob\P\}$. Suppose the objects $x$ occurred in this expression form a subset $\S\subset\Ob\P$. In this way $Z$ uniquely determines $\S$. On the other hand, given $\S\subset\Ob\P$, we may always define
$$
Z(\S)=\{\p^x \bigm{|}x\in\S\subset\Ob\P\}\subset\Spc\K(\P).
$$
The fourth correspondence can be constructed by composing the other three. To be more precise, we may associate to $Z(\S)$ a tensor ideal $\I=\langle S_x\bigm{|} x\not\in\S\rangle$. Conversely, let $\I$ be a tensor ideal of $\K(\P)$. We assign to it a subset of $\Spc\K(\P)$ by
$$
\I \mapsto \{\p^x\bigm{|}S_x\hotimes\I=0\}.
$$
The latter is identified with $\supp(\I)$.
\end{example}

What we conclude from our discussions is that, although the finite generation of the cohomology ring $\H^*(\C)$ is highly unpredictable, one always has a Noetherian space $\Spc\K(\C)$ which reveals plenty of representation-theoretic information of $k\C$, whenever $\C$ is finite EI.

\section{Transporter category algebras}

One of the distinguished feature of a finite transporter category algebra is that it is Gorenstein, so we may consider a different kind of tensor triangulated categories, $\CM k(G\propto\P)$.

Given a finite transporter category algebra $k(G\propto\P)$, Theorem 3.3.1 reads as follows
$$
\Spc\K(G\propto\P)\simeq\biguplus_{[x]\subset\Ob(G\propto\P)}\Spc\K(G\propto[x])\simeq\biguplus_{[x]\subset\Ob(G\propto\P)}\Spec^h\H^{\cdot}(G_x),
$$
where $G\propto[x]=(G\propto\P)_x$ and $G_x=\Aut_{G\propto\P}(x)$.

\subsection{Gorenstein property}

Suppose $\P$ is a finite $G$-poset. Then we have a category algebra, the incidence algebra, $k\P$, as well as a category algebra $k(G\propto\P)$. The algebra $k\P$ inherits an obvious $G$-action and thus we can construct a skew group algebra $k\P[G]$, an algebraic semi-direct product between $kG$ and $kP$. We have the following result from \cite{X2}.

\begin{lemma} Let $\P$ be a finite $G$-poset for a finite group $G$. Then
$$
k(G\propto\P)\cong k\P[G]
$$
is Gorenstein.
\end{lemma}

This identification allows us to consider the maximal Cohen-Macaulay modules over $k(G\propto\P)$, see \cite{B, H} for the explicit definition. Indeed the tensor product on $\K(G\propto\P)=\D^b(k(G\propto\P)\mbox{-}\mod)$ passes to the quotient category
$$
\CM k(G\propto\P)=\D^b(k(G\propto\P)\mbox{-}\mod)/\D^b(k(G\propto\P)\mbox{-}\proj),
$$
for $\D^b(k(G\propto\P)\mbox{-}\proj)$ is a tensor ideal by the following lemma of \cite{X1}.

\begin{lemma} $\M\in k(G\propto\P)$-$\mod$ is of finite projective (equivalently finite injective) dimension if and only if $\M_x\in k(G\propto[x])$-$\mod$ is projective (equivalently injective) for all $[x]\subset\Ob(G\propto\P)$.
\end{lemma}

In the previous lemma, $\M_x\in k(G\propto[x])$-$\mod$ being projective is equivalent to $\M(x)\in kG_x$-$\mod$ being projective.

\subsection{Stable categories of maximal Cohen-Macaulay modules} We continue to compute the spectrum of the tensor triangulated category $\CM k(G\propto\P)$ here. By Lemma 4.1.2, if all automorphism group algebras are semi-simple, the category $\CM k(G\propto\P)$ is $\{0\}$.

Let us first settle some terminology. Let $\Q\subset\P$ be a $G$-subposet. Then $G\propto\Q$ becomes a full subcategory of $G\propto\P$. It is easy to see that every full subcategory of $G\propto\P$ is a transporter category on some $G$-subposet $\Q$. In light of this, we shall simplify our earlier notations $\R_{G\propto\P}$ and $\Inc_{G\propto\Q}$ to $\R_{\Q}$ and $\Inc_{\Q}$, respectively. Particularly $\R_{[x]}=\R_{G\propto[x]}$ and $\Inc_{[x]}=\Inc_{G\propto[x]}$.

\begin{theorem} There is a homeomorphism
$$
\Spc\CM k(G\propto\P)\simeq\biguplus_{[x]\subset\Ob(G\propto\P)}\Proj\H^{\cdot}(G_x).
$$
It results in a bijection between the set of tensor ideals in $\CM k(G\propto\P)$ and the set of specialization closed subsets of $\biguplus_{[x]\subset\Ob(G\propto\P)}\Proj\H^{\cdot}(G_x)$.

\begin{proof} Consider the following localizing sequence
$$
\D^b(k(G\propto\P)\mbox{-}\proj) \to \K(G\propto\P) \to \CM k(G\propto\P).
$$
From \cite[Proposition 3.11]{Ba1} we know that $\Spc\CM k(G\propto\P)$ is homeomorphic to a subspace
$$
X=\{\p\in\Spc\K(G\propto\P)\bigm{|}\D^b(k(G\propto\P)\mbox{-}\proj)\subset\p\}\subset\Spc\K(G\propto\P).
$$
Since $\R_{[x]}\D^b(k(G\propto\P)\mbox{-}\proj)=\D^b(k(G\propto[x])\mbox{-}\proj)\subset\R_{[x]}\p$, the only $\p\in\Spc\K(G\propto\P)$, not belonging to this subspace, must satisfy $\R_{[x]}\p=0$ for the unique $[x]$ due to the reason that $\D^b(k(G\propto[x])\mbox{-}\proj)$ is the smallest non-zero tensor ideal of $\K(G\propto[x])=\D^b(k(G\propto[x])\mbox{-}\mod)$. These primes ideals are indexed by the set of isomorphism classes of objects $\{[x] \bigm{|} x\in\Ob(G\propto\P)\}$. Let us denote
$$
Y=\Spc\K(G\propto\P)\backslash X.
$$
Analogous to \cite[Proposition 8.5]{Ba2} we have a commutative diagram
$$
\xymatrix{\Spc\CM k(G\propto\P) \ar@{^{(}->}[r] \ar[d]^{\simeq}_{\widetilde{\R}} & \Spc\K(G\propto\P) \ar[d]^{\simeq}_{\widetilde{\R}} & Y \ar@{_{(}->}[l] \ar[d]^{\simeq}_{\widetilde{\R}}\\
\biguplus_{[x]}\Spc\CM k(G\propto[x]) \ar@{^{(}->}[r] & \biguplus_{[x]}\Spc\K(G\propto[x]) \ar[d]_{\biguplus\rho^{\cdot}}^{\simeq} & \biguplus_{[x]}\{\H^+(G\propto[x])\} \ar@{_{(}->}[l] \ar[d]_{\biguplus\rho^{\cdot}}^{\simeq}\\
\biguplus_{[x]}\Proj\H^{\cdot}(G_x) \ar@{^{(}->}[r] \biguplus_{[x]}\ar[u]^{\biguplus\varphi}_{\simeq} & \biguplus_{[x]}\Spec^h\H^{\cdot}(G_x) & \biguplus_{[x]}\{\H^+(G_x)\} \ar@{_{(}->}[l]}
$$
In summary we have the following homeomorphisms
$$
\Spc\CM k(G\propto\P) \simeq \biguplus_{[x]\subset\Ob(G\propto\P)}\Spc\CM k(G\propto[x]) \simeq \biguplus_{[x]\subset\Ob(G\propto\P)}\Proj\H^{\cdot}(G_x).
$$

\end{proof}
\end{theorem}

Similar to remarks after Theorem 3.3.1, the structure presheaves on the two topological spaces $\Spc\CM k(G\propto\P)$ and $\biguplus_{[x]\subset\Ob(G\propto\P)}\Proj\H^{\cdot}(G_x)$ are quite different. For instance, the global section on the space $\Spc\CM k(G\propto\P)$ is the Tate cohomology ring $\widehat{\H}^*_G(B\P, k)\cong\widehat{\Ext}^*_{k(G\propto\P)}(\k,\k)$, of the transporter category $G\propto\P$. By contrast the global section on $\biguplus_{[x]\subset\Ob(G\propto\P)}\Proj\H^{\cdot}(G_x)$ is the direct product of Tate cohomology rings of all $G_x$, $[x]\subset\Ob(G\propto\P)$.

\end{document}